\documentclass[11pt]{article}

\usepackage{amsmath}
\usepackage{amsthm}
\usepackage{amssymb}
\usepackage{mathptmx}
\usepackage{concrete}
\def\pf{{\bf Proof. }}

\def\bel#1{\begin{equation}\label #1}

\def\l{\lambda}

\def\D{\Delta}

\def\nd{\noindent}
\def\p{\partial}

\def\BB{{\Bbb B}}
\def\hf{\hfill{$\Box$}}
\def\<{\leq}
\def\>{\geq}

\begin{document}

\title{\bf Gradient estimate of a Dirichlet eigenfunction on a compact manifold with boundary}
\author{Yiqian Shi$^{*}$ and Bin Xu$^{*\dagger}$}
\date{}
\maketitle

\nd{\small {\bf Abstract.}\ \ Let $e_\l(x)$ be an eigenfunction
with  respect to the  Dirichlet Laplacian $\Delta_N$ on a compact
Riemannian manifold $N$ with boundary: $\Delta_N e_\l=\l^2 e_\l$
in the interior of $N$ and $e_\l=0$ on the boundary of $N$.  We
show the following gradient estimate of $e_\l$: for every $\l\geq
1$, there holds $\l\|e_\l\|_\infty/C\leq \|\nabla
e_\l\|_\infty\leq C\l \|e_\l\|_\infty$, where $C$ is a positive
constant depending only on $N$. In the proof, we use a basic
geometrical property of nodal sets of eigenfunctions and elliptic
apriori estimates. }

\footnote{\hspace{-0.7cm} $^*$ Department of Mathematics,
University of Science
and Technology of China, Hefei 230026 China.\\
$^\dagger$ E-mail of the correspondent author: bxu@ustc.edu.cn}

\nd {\small {\bf Mathematics Subject Classification (2000):}\
Primary 35P20; Secondary 35J05}

\nd {\small {\bf Key Words:}\ Dirichlet Laplacian,  eigenfunction,
gradient estimate}

\section{Introduction}

 \quad  Let $(N,\ g)$ be an $n$-dimensional
compact smooth Riemannian manifold with smooth boundary $\p N$ and
$\Delta_N$ the positive Dirichlet Laplacian on $N$. Let $L^2(N)$
be the space of square integrable functions on $N$ with respect to
the Riemannian density $dv(N)=\sqrt{{\bf g}(x)}\,dx:= \sqrt{\det\
(g_{ij})}\ dx$. Let $e_1(x),\,e_2(x),\,\cdots$ be a complete
orthonormal basis in $L^2(N)$ for the Dirichlet eigenfunctions of
$\D_N$ such that $0<\l_1^2< \l_2^2\leq  \l_3^2\leq\,\cdots$ for
the corresponding eigenvalues, where $e_j(x)$ ($j=1,2,\dots$) are
real valued smooth function on $N$ and $\l_j$ are positive
numbers. Also, let ${\bf e}_j$ denote the projection of $L^2(N)$
onto the 1-dimensional space ${\bf C}e_j$. Thus , an $L^2$
function $f$ can be written as $f=\sum_{j=0}^\infty {\bf e}_j(f)$,
where the partial sum converges in the $L^2$ norm. Let $\l$ be a
positive real number $\geq  1$. We define the spectral function
and the unit band spectral projection operator $\chi_\l$ as
follows:
\[e(x,y,\l):=\sum_{\l_j\leq \l}\,e_j(x)\,e_j(y),\]
\[{\chi}_\l f:=\sum_{\l_j\in (\l,\,\l+1]}{\bf e}_j(f)\ .\]

Grieser \cite{Gr} and Sogge \cite{S3} proved the $L^\infty$
estimate of $\chi_\l$,
\begin{equation}
\label{equ:sup}
 ||{ \chi}_\l f||_\infty\leq C\l^{(n-1)/2} ||f||_2\,
\end{equation}
where $||f||_r$ ($1\leq r\leq \infty$) means the $L^r$ norm of the
function $f$ on $N$. In the whole of this paper $C$ denotes a
positive constant which depends only on $N$ and may take different
values at different places, if there is no otherwise stated. The
idea of Grieser and Sogge is to use the standard wave kernel
method outside a boundary layer of width $C \lambda^{-1}$ and a
maximum principle argument inside that layer. By using the maximum
principle argument and the estimate (\ref{equ:sup}), Xu \cite{X2}
proved the gradient estimate of $\chi_\lambda$
\begin{equation}
\label{equ:2infgrad} ||\nabla\,  { \chi}_\l f||_\infty\leq
C\l^{(n+1)/2} ||f||_2\ .
\end{equation}
Here $\nabla$ is the Levi-Civita connection on $N$. In particular,
$\nabla f=\sum_j\, g^{ij}\partial f/\partial x_j$ is the gradient
vector field of a $C^1$ function $f$, the square of whose length
equals $\sum_{i,j} g^{ij}(\p f/\p x_i)(\p f/\p x_j)$. One of his
motivation is to prove the H{\" o}rmander multiplier theorem on
compact manifolds with boundary.  Seeger and Sogge \cite{SS}
firstly proved that theorem by using the parametrix of the wave
kernel on manifolds without boundary. All the results mentioned in
the introduction have their analog on compact manifolds without
boundary. See the details in the introduction of Shi-Xu \cite{SX}
and the references therein.  In general, the method used in
manifolds without boundary is not valid for the problems on
manifolds with boundary. In particular, on manifolds with boundary
the H{\" o}rmander multiplier theorem cannot be obtained by the
standard pseudo-differential operator calculus as done on
manifolds without boundary, since the square root of the Dirichlet
Laplaican is not a pseudo-differential operator any more and one
cannot obtain the $L^\infty$ bounds for $\chi_\l$ and $\nabla
\chi_\l$ only by using the Hadamard parametrix of the wave kernel.

In the paper, by rescaling $\chi_\l f$ at the scale of $\l^{-1}$
both outside and inside the boundary layer of width $C\l^{-1}$, we
obtain by elliptic apriori estimates a slightly
stronger estimate than (\ref{equ:2infgrad}) as follows:\\

\nd {\bf Theorem 1.1.} {\it Let $f$ be a square integrable
function on $N$. Then, for every $\l\geq 1$, there holds
\begin{equation}
\label{equ:grad2}
 \|\nabla \chi_\l f\|_\infty\leq C\left(\l \|\chi_\l
 f\|_\infty+\l^{-1}\|\Delta_N\,\chi_\l f\|_\infty\right).
\end{equation}}

\nd {\bf Remark 1.1.} Putting $f(\cdot)=\sum_{\l_j\in
(\l,\,\l+1]}\, e_j(x)e_j(\cdot)$ in (\ref{equ:sup}), we obtain the
uniform estimate of eigenfunctions for all $x\in N$,
\[\sum_{\l_j\in
(\l,\,\l+1]}\, |e_j(x)|^2\leq C\l^{n-1}.\] Actually it can be
proved that those two estimates are equivalent. By using the
Cauchy-Scahwarz inequality, we obtain the gradient estimate
(\ref{equ:2infgrad}) from Theorem 1 together with the above
inequality. Similarly, that estimate is equivalent to the uniform
estimate for all $x\in N$,
\begin{equation}
\label{equ:gradpointwise}
 \sum_{\l_j\in (\l,\,\l+1]}\, |\nabla
e_j(x)|^2\leq C\l^{n-1}.\end{equation}

\nd{\bf Remark 1.2.}  By the finite propagation speed of the wave
equation, the asymptotic formula of derivatives of the spectral
function $e(x,y,\l)$ in Theorem 1 \cite{XuB}  which is proved by
the standard wave kernel method, also holds for each interior
point $x$ of $N$. Much more general asymptotic formulae are given
in Theorems 1.8.5 and 1.8.7 of Safarov-Vassiliev \cite{SV}. In
particular, we have the following asymptotic formula that as
$\l\to\infty$
\[\sum_{\l_j\leq
\l}|\nabla\,e_j(x)|^2=\frac{n\,\l^{n+2}}{2\,(4\pi)^{n/2}\,\Gamma(2+\frac{n}{2})}+{\rm
O}_x(\l^{n+1}),\] where the constant in the reminder term ${\rm
O}_x(\l^{n+1})$ depends on the distance of $x$ to the boundary of
$N$. Hence, the exponents of $\l$ in estimates
(\ref{equ:2infgrad}) and (\ref{equ:gradpointwise}) are sharp at
$x$ as $\l\to\infty$. For each point $z$ on the boundary of $N$,
Ozawa \cite{Oz} used the heat kernel method to show the asymptotic
formula that as $\l\to\infty$,
\[\sum_{\l_j\leq
\l}\left|\frac{\p e_j}{\p
\nu}(z)\right|^2=\frac{\l^{n+2}}{(4\pi)^{n/2}\,\Gamma(2+\frac{n}{2})}+{\rm
o}(\l^{n+2}),\] where $\nu$ is the unit outer normal vector field
on the boundary of $N$. Hence, (\ref{equ:2infgrad}) and
(\ref{equ:gradpointwise}) are also sharp on the boundary. \\

\nd {\bf Corollary 1.1.} {\it Let $e_\l(x)$ be an eigenfunction
with respect to the positive Dirichlet Laplacian $\Delta_N$ on
$N$: $\Delta_M e_\l=\l^2 e_\l$ in the interior of $N$ and $e_\l=0$
on the boundary of $N$. Then, for every $\l\geq 1$, there holds
the upper bound estimate of $\nabla e_\l$ {\rm :}
\[\|\nabla e_\l\|_\infty\leq C\l\,\|e_\l\|_\infty.\]}

\nd \pf Putting $f=e_\l$ in the estimate (\ref{equ:grad2}), we
obtain the
corollary. \hf \\

Actually, by the basic geometric property of nodal sets of an
eigenfunction, we can find a complete picture for the $L^\infty$
norm of $\nabla e_\l$ in the following:\\

\nd {\bf Theorem 1.2.} {\it Let $e_\l(x)$ be an eigenfunction with
respect to the Dirichlet Laplacian operator $\Delta_N$ on $N$
without boundary: $\Delta_N e_\l=\l^2 e_\l$. Then, for every
$\l\geq 1$, there holds
\begin{equation}
\label{equ:grad} \l\|e_\l\|_\infty/C\leq \|\nabla
e_\l\|_\infty\leq C\l \|e_\l\|_\infty.
\end{equation}
}\\

\nd {\bf Remark 1.3.} The authors \cite{SX} proved the analog of
Theorems 1 and 2 on compact manifolds without boundary. The proof
in this paper is more complicated than there because we need to do
analysis at points near the boundary. We believe that there also
hold the analog for $k$-covariant derivatives $\nabla^k\,\chi_\l
f$ and $\nabla^k \,e_\l$ on $N$. We plan to
discuss this question in a future paper.\\

\nd {\bf Remark 1.4.} Let $\psi_j$ be the normal derivative of
$e_j$ at the boundary $\p N$ of $N$. The lower bound estimate
\[\|\psi_j\|_{L^\infty(\p N)}\geq C\|e_j\|_{L^\infty(N)}\]
does not hold in general. Using Examples 3-5 in Hassell-Tao
\cite{HT} and doing a little bit more computations, we can see
that the above estimate does not hold on the flat cylinder, the
hemisphere and the spherical cylinder. We hope to find a
sufficient condition
for the lower bound estimate in a future work.\\

We conclude the introduction by explaining the organization of
this paper. In Section 2, we show the lower bound of the gradient
$\nabla\, e_\l$ by the basic geometrical property of the nodal set
of eigenfunctions. In Section 3 we use the rescaling method and
the H{\" o}lder estimate about elliptic PDEs to show
(\ref{equ:grad2}) and the upper bound part of (\ref{equ:grad}).
The point is to do the rescaling both outside and inside the
boundary layer of width $C\,\l^{-1}$.

\section{Lower bound of $\nabla e_\l$}

\quad \ \ The nodal set of an eigenfunction $e_\l$ of $\D_N$ is
the zero set
\[Z_{e_\l}:=\{x\in N: e_\l(x)=0\}.\]
A connected component of the open set $N\backslash Z_{e_\l}$ is
called a nodal domain of the eigenfunction $e_\l$. We have the
same definition for manifolds without boundary.\\

\nd {\bf Lemma 2.1.} (Br{\" u}ning \cite{Br}) {\it Let $M$ be a
compact Riemannian manifold without boundary. Let $\l\geq 1$ and
$e_\l$ be an eigenfunction of the positive Laplacican $\D_M$:
$\D_M\,e_\l=\l^2 e_\l$. Then there exists a constant $C$ only
depending on $M$ such that each geodesic ball of radius $C/\l$ in
$M$ must intersect the nodal set $Z_{e_\l}$ of $e_\l$.}

\nd\pf A proof written in English is given by Zelditch in pp. 579-580 of \cite{Z}.\hf\\

We need a manifold-with-boundary version of Lemma 1 as follows:\\

\nd {\bf Lemma 2.2.}  {\it Let $\l>0$ and $\l^2$ be greater than
the smallest eigenvalue $\l_1^2$ of the Dirichlet Laplacian
$\D_N$. Let $e_\l$ be an eigenfunction of $\D_N$: $\D_N e_\l=\l^2
e_\l$ in the interior ${\rm Int}(N)$ of $N$ and $e_\l=0$ on the
boundary $\p N$ of $N$. Then there exists a positive constant $D$
only depending on $N$ such that each geodesic ball of radius
$D/\l$ contained in ${\rm Int}(N)$ must
intersect the nodal set $Z_{e_\l}$ of $e_\l$.}\\

\nd \pf We here adapt the proof of Zelditch \cite{Z} with a slight
modification.\\

{\it Step 1}\quad We show the following fact: There exists a
constant $C$ such that for each interior point $p$ of $N$ and each
positive number $r>0$ satisfying that the distance $d(p,\,\p N)<r$
from $p$ to $\p N$ is less than $r$, the smallest Dirichlet
eigenvalue $\l_1^2\bigl(B(p,\,r)\bigr)$ of the geodesic ball
$B(x,\,r)$ is bounded from above by $C/r^2$. Since $d(p,\,\p
N)>r$, we may assume that there exists a geodesic normal
coordinate chart $(x_1,\dots, x_n)$ on the ball $B(p,\,r)$. Let
$g$ be the Riemannian metric of $N$ on $B(p,\,r)$ with
coefficients $g_{ij}=g(\p/\p x_i,\,\p/\p x_j)$ and $g_0$ be the
Euclidean metric $dx_1^2+\cdots+dx_n^2$ on $B(x,\,r)$. Take
$0<c_1<1$ depending only on $N$ so that the Euclidean ball
$B(p,c_1r\,;\,g_0)$ is contained in the metric ball $B(p,r\,;\,g)$
of $N$. Since $c_1<1$, by the definition of Rayleigh quotient,
\[\l_1^2\bigl(B(p,r\,;\,g)\bigr)\leq
\l_1^2\bigl(B(p,c_1r\,;\,g)\bigr).\] Since $N$ is compact, by
comparing Rayleigh quotients, there exists $c_2>0$ depending only
on $N$ such that
\[\l_1^2\bigl(B(p,c_1r\,;\,g)\bigr)\leq c_2
\l_1^2\bigl(B(p,c_1r\,;\,g_0)\bigr).\] On the other hand, by
change of variables, we have
\[\l_1^2\bigl(B(p,c_1r\,;\,g_0)\bigr)=\frac{\l_1^2\bigl(B(p,1\,;\,g_0)\bigr)}{(c_1r)^2}.\]
Combining the above three inequalities and setting
$C=c_2c_1^{-2}\l_1^2\bigl(B(p,1\,;\,g_0)\bigr)$, we complete the
proof.\\

{\it Step 2} Take a geodesic ball $B(p,\,r)$ in $N$ such that
$d(p,\,\p N)>r$. Suppose that it is disjoint from the nodal set
$Z_{e_\l}$. Then it is completely contained in a nodal domain
$D_j$ of $e_\l$. But $\l^2=\l_1^2(D_j)\leq
\l_1^2\bigl(B(p,\,r)\bigr)\leq
C/r^2$. Hence, $r\leq \sqrt{C}{\l}$. Taking $D=2\sqrt{C}$, we complete the proof.\hf  \\

\nd {\sc Proof of the lower bound part of Theorem 2}\quad Take a
point $x$ in $N$ such that $|e_\l(x)|=\|e_\l\|_\infty$. By the
Dirichlet boundary condition, the distance $d$ from $x$ to $\p N$
is positive.

{\it Case 1} \quad Assume $d>D/\l$. Then there exists point $y$ in
the geodesic ball $B(x,\,D/\l)$ with center $x$ and radius $D/\l$
such that $e_\l(y)=0$. We may assume $\l$ so large that there
exists a geodesic normal chart $(r, \theta)\in [0,\,D/\l]\times
{\Bbb S}^{n-1}(1)$ in the ball $B(x,\,D/\l)$. By the mean value
theorem, there exists a point $z$ on the geodesic segment
connecting $x$ and $y$ such that
\[\left|\frac{\p e_\l}{\p r}(z)\right|\geq
\frac{\l}{D}|e_\l(x)|=\frac{\l}{D}\|e_\l\|_\infty.\]

{\it Case 2} \quad Assume $d\leq D/\l$. We may assume $\l$ so
large that there exists a unique geodesic $\gamma:[0,\,d]\to N$ of
arc length parameter connecting $x$ and $\p N$,
\[\gamma(0)=x,\quad \gamma(d)\in \p N.\]
Since $e_\l(\gamma(d))=0$, by the mean value theorem, there exists
$t_0$ in $(0,\,d)$ such that
\[\left|\frac{d
e_\l\bigl(\gamma(t)\bigr)}{dt}(t_0)\right|=\frac{|e_\l(x)|}{d}\geq
\frac{\l}{D}\,\|e_\l\|_\infty.\] \hf

\section{Estimate for $\nabla \chi_\l f$}

\subsection{Outside the boundary layer}

\quad\ \ Recall the principle: {\it On a small scale comparable to
the wavelength $1/\l$, the eigenfunction $e_\l$ behaves like a
harmonic function.} It was developed in H. Donnelly  and C.
Fefferman \cite{DF1} \cite{DF2} and N. S. Nadirashvili \cite{Na}
and was used extensively there. Recently Mangoubi \cite{Ma1}
applied this principle to studying the geometry of nodal domains
of eigenfunctions. In this section, for a square integrable
function $f$ on $N$ we give a modification of this principle,
which can be applied to the Poisson equation
\[\D_N\,\chi_\l f=\sum_{\l_j\in (\l,\,\l+1]}\, \l_j^2 {\bf
e}_j(f)\quad {\rm in}\quad {\rm Int}(N)\] with the Dirichlet
boundary condition $\chi_\l f=0$ on $\p N$. In particular, in this
subsection, we do the analysis outside the boundary layer
$L_{1/\l}=\{z\in N:\,d(z,\,\p N)\leq 1/\l\}$ of width $1/\l$.

Take an arbitrary point $p$ with $d(p,\,\p N)\geq 1/\l$. We may
assume that $1/\l$ is sufficiently small such that there exists a
geodesic normal coordinate chart $(x_1,\dots, x_n)$ on the
geodesic ball $B(p,\,2/\l)$ in $N$. In this chart, we may identify
the ball $B(p,\,2/\l)$ with the $n$-dimensional Euclidean ball
${\Bbb B}(2/\l)$ centered at the origin $0$, and think of the
function $\chi_\l f$ in $B(p,\,2/\l)$ as a function in ${\Bbb
B}(2/\l)$.  Our aim in this subsection is to show the inequality
\begin{equation}
\label{equ:outside2}
 |(\nabla \chi_\l f)(p)|\leq
C\left(\l\|\chi_\l f\|_{L^\infty\bigl({\Bbb
B}(2/\l)\bigr)}+\l^{-1}\|\D_N\,\chi_\l f\|_{L^\infty\bigl({\Bbb
B}(2/\l)\bigr)}\right).
\end{equation}
For simplicity of notions, we rewrite $u=\chi_\l f$ and
$v=\D_N\chi_\l f$ in what follows. The Poisson equation satisfied
by $u$ in ${\Bbb B}(1/\l)$ can be written as
\[-\frac{1}{\sqrt{g}}\sum_{i,j}\,\p_{x_i}\left(g^{ij}\sqrt{g}\p_{x_j}
u\right)=v.\]

Consider the rescaled functions $u_\l(y)=u(y/\l)$ and
$v_\l(y)=v(y/\l)$ in the ball ${\Bbb B}(2)$. The above estimate we
want to prove is equivalent to its rescaled version
\begin{equation}
\label{equ:scale}
 |(\nabla u_\l)(0)|\leq
C\left(\|u_\l\|_{L^\infty\bigl({\Bbb
B}(2)\bigr)}+\l^{-2}\|v_\l\|_{L^\infty\bigl({\Bbb
B}(2)\bigr)}\right).
\end{equation}
On the other hand, the rescaled version of the Poisson equation
has the expression,
\begin{equation}
\label{equ:PoissonScaled}
\sum_{i,j}\,\p_{y_i}\left(g^{ij}_\l\sqrt{g_\l}\p_{y_j}
u_\l\right)=-\l^{-2}\,\sqrt{g_\l}\,v_\l,\end{equation} where
$g_{ij,\,\l}(y)=g_{ij}(y/\l)$, $g^{ij}_\l(y)=g^{ij}(y/\l)$ and
$\sqrt{g_r}(y)=(\sqrt{g})(y/\l)$.

For each $0<\alpha<1$, there exists $K>0$ such that the $C^\alpha$
norm of the coefficients $g^{ij}_\l\sqrt{g_\l}$, $\sqrt{g_\l}$ in
$\BB(2)$ are bounded uniformly from above by $K$, and the smallest
eigenvalue of the $n\times n$ matrix $(g^{ij}_\l\sqrt{g_\l})_{ij}$
in $\BB(2)$ bounds from below by $1/K$, for all $\l\geq 1$. By
Theorem 8.32 in page 210 of Gilbarg-Trudinger \cite{GT}, there
exists constant $C=C(n, \, \alpha,\, K)$ such that
\[\|u_\l\|_{C^{1,\,\alpha}\bigl({\Bbb B}(1)\bigr)}\leq
C\left(
\|u_\l\|_{L^\infty\bigl(\BB(2)\bigr)}+\l^{-2}\,\|v_\l\|_{L^\infty\bigl(\BB(2)\bigr)}\right),
\]
This is stronger than the estimate (\ref{equ:scale}). Therefore,
we complete the proof of Theorem 1.1 outside the boundary layer
$L_\l$.

\subsection{Inside the boundary layer}

\quad\ \ Using the notions in subsection 3.1, We are going to
prove the following estimate:
\begin{equation}
\label{equ:inside} \|\nabla u\|_{L^\infty(L_{1/\l})}\leq C\left(\l
\|u\|_\infty+\l^{-1}\|v\|_\infty\right),
\end{equation}
with which combining (\ref{equ:outside2}) completes the proof of
Theorem 1.1.

We may assume that $\l$ is sufficiently large so that there exists
a geodesic normal coordinate chart $(z',\,z_n)$ on the boundary
layer $L_{3/\l}=\{p\in N:\,d(p,\,\p N)\leq 3/\l\}$ with respect to
the boundary $\p N$. Hence, for each point $(z',\,z_n)\in
L_{3/\l}$, we have $0\leq z_n\leq 3/\l$ and
\[d\bigl((z',\,z_n),\,\p N\bigr)=z_n.\]
For each point $q\in \p N$ and $r>0$, denote by $B_+(q,\,r)$ the
set of points of $N$ with distance less than $r$ to $q$. Denote by
$\BB_+(r)$ the upper half Euclidean ball
\[\{x=(x_1,\dots,x_n)\in {\bf R}^n:|x|<r,\,x_n\geq 0\}\]
centered at the origin and with radius $r$. Then, for each $q\in
\p N$, there exists a geodesic normal chart on $B_+(q,\,3/\l)$
such that the exponential map $\exp_q$ at $q$ gives a
diffeomorphism from $\BB_+(3/\l)$ onto $B_+(q,\,3/\l)$.

Since $\{B_+(q,\,2/\l):q\in \p N\}$ forms an open cover of
$L_{1/\l}$, the question can be reduced to showing the analog of
(\ref{equ:inside}) on $B_+(q,\,2/\l)$ for each $q$. We only need
to prove its equivalent rescaled version,
\begin{equation}
\label{equ:insidescaled} \|\nabla
u_\l\|_{L^\infty\bigl(\BB_+(2)\bigr)}\leq C\left(
\|u_\l\|_{L^\infty\bigl(\BB_+(3)\bigr)}+\l^{-2}\|v_\l\|_{L^\infty\bigl(\BB_+(3)\bigr)}\right),
\end{equation}
where $u_\l$ and $v_\l$ are the the rescaling function of $u$ an
$v$, respectively. Observe that $u_\l$ and $v_\l$ satisfy the
Poisson equation $(\ref{equ:PoissonScaled})$ in the upper half
Euclidean ball $\BB(3)$ and the Dirichlet boundary condition,
\[u_\l=0\quad {\rm on\ the\ portion}\quad \{x\in \BB(3):x_n=0\}\]
of the boundary of $\BB_+(3)$.  For each $0<\alpha<1$, there
exists $K>0$ such that the $C^\alpha$ norm  of the coefficients
$g^{ij}_\l\sqrt{g_\l}$, $\sqrt{g_\l}$ in $\BB_+(3)$ are bounded
uniformly from above by $K$, and the smallest eigenvalue of the
$n\times n$ matrix $(g^{ij}_\l\sqrt{g_\l})_{ij}$ in $\BB_+(3)$
bounds from below by $1/K$, for all $\l\geq 1$. By Theorem 8.36 in
page 212 of Gilbarg-Trudinger \cite{GT}, there exists constant
$C=C(n, \, \alpha,\, K)$ such that
\[\|u_\l\|_{C^{1,\,\alpha}\bigl({\Bbb B}_+(2)\bigr)}\leq
C\left(
\|u_\l\|_{L^\infty\bigl(\BB_+(3)\bigr)}+\l^{-2}\,\|v_\l\|_{L^\infty\bigl(\BB_+(3)\bigr)}\right).
\]
This is a stronger estimate than (\ref{equ:insidescaled}). \\

\nd  {\bf Acknowledgements}\quad Yiqian Shi is supported in part
by the National Natural Science Foundation of China (No. 10671096
and No. 10971104) and Bin Xu by the National Natural Science
Foundation of China (No. 10601053 and No. 10871184).

\end{document}